\documentclass[12pt]{amsart}
\usepackage{amsfonts,amssymb,color}
\usepackage[mathscr]{eucal}
\usepackage{amsmath, amsthm}
\usepackage{mathrsfs}
\usepackage{amsbsy}
\usepackage{dsfont}
\usepackage{bbm}
\makeindex
\input xypic
\xyoption{all}
\begin{document}
\baselineskip = 5mm
\newcommand \lra {\longrightarrow}
\newcommand \hra {\hookrightarrow}
\newcommand \ZZ {{\mathbb Z}} 
\newcommand \NN {{\mathbb N}} 
\newcommand \QQ {{\mathbb Q}} 
\newcommand \RR {{\mathbb R}} 
\newcommand \CC {{\mathbb C}} 
\newcommand \bcA {{\mathscr A}}
\newcommand \bcB {{\mathscr B}}
\newcommand \bcC {{\mathscr C}}
\newcommand \bcD {{\mathscr D}}
\newcommand \bcE {{\mathscr E}}
\newcommand \bcF {{\mathscr F}}
\newcommand \bcI {{\mathscr I}}
\newcommand \bcJ {{\mathscr J}}
\newcommand \bcM {{\mathscr M}}
\newcommand \bcP {{\mathscr P}}
\newcommand \bcS {{\mathscr S}}
\newcommand \bcT {{\mathscr T}}
\newcommand \bcU {{\mathscr U}}
\newcommand \bcX {{\mathscr X}}
\newcommand \bcY {{\mathscr Y}}
\newcommand \bcZ {{\mathscr Z}}
\newcommand \C {{\mathscr C}}
\newcommand \im {{\rm im}}
\newcommand \Hom {{\rm Hom}}
\newcommand \colim {{{\rm colim}\, }} 
\newcommand \End {{\rm {End}}}
\newcommand \Aut {\rm Aut}
\newcommand \coker {{\rm {coker}}}
\newcommand \id {{\rm {id}}}
\newcommand \supp {{\rm {Supp}}\, }
\newcommand \CHM {{\bcC \! \bcM }}
\newcommand \SP {{\bcS \! \bcP }}
\newcommand \DM {{\mathscr D\! \mathscr M}}
\newcommand \MM {{\mathscr M\! \mathscr M}}
\newcommand \uno {{\mathbbm 1}}
\newcommand \Le {{\mathbbm L}}
\newcommand \ptr {{\pi _2^{\rm tr}}}
\newcommand \Ob {{\rm Ob}}
\newcommand \PR {{\mathbb P}} 
\newcommand \AF {{\mathbb A}} 
\newcommand \Spec {{\rm {Spec}}}
\newcommand \Pic {{\rm {Pic}}}
\newcommand \Alb {{\rm {Alb}}}
\newcommand \Jac {{\rm {Jac}}}
\newcommand \Sym {{\rm {Sym}}}
\newcommand \Corr {{CH}}
\newcommand \cha {{\rm {char}}}
\newcommand \tr {{\rm {tr}}}
\newcommand \trdeg {{\rm {tr.deg}}}
\newcommand \gm {{\mathfrak {m}}}
\newcommand \gp {{\mathfrak {p}}}
\def\blue {\color{blue}}
\def\red{\color{red}}
\newtheorem{theorem}{Theorem}
\newtheorem{lemma}[theorem]{Lemma}
\newtheorem{sublemma}[theorem]{Sublemma}
\newtheorem{corollary}[theorem]{Corollary}
\newtheorem{example}[theorem]{Example}
\newtheorem{exercise}[theorem]{Exersize}
\newtheorem{proposition}[theorem]{Proposition}
\newtheorem{remark}[theorem]{Remark}
\newtheorem{notation}[theorem]{Notation}
\newtheorem{definition}[theorem]{Definition}
\newtheorem{conjecture}[theorem]{Conjecture}
\newtheorem{claim}[theorem]{Claim}
\newenvironment{pf}{\par\noindent{\em Proof}.}{\hfill\framebox(6,6)
\par\medskip}
\title[Transcendence degree of zero-cycles]
{\bf Transcendence degree of zero-cycles and the structure of Chow motives}
\author{S. Gorchinskiy, V. Guletski\u \i }


\begin{abstract}
\noindent In the paper we introduce a transcendence degree of a zero-cycle on a smooth projective variety $X$ and relate it to the structure of the motive of $X$. In particular, we show that in order to prove Bloch's conjecture for a smooth projective complex surface $X$ of general type with $p_g=0$ it suffices to prove that one single point of a transcendence degree $2$ in $X(\CC )$, over the minimal subfield of definition $k\subset \CC $ of $X$, is rationally equivalent to another single point of a transcendence degree zero over $k$. This can be of particular interest in the context of Bloch's conjecture for those surfaces which admit a concrete presentation, such as Mumford's fake surface, see \cite{Mumford}.
\end{abstract}





\maketitle

\section{Introduction}
\label{s-intro}

Since \cite{Bloch} we know that the generic point, considered as a zero-cycle, plays an important role in the study of algebraic cycles on a smooth projective variety $X$ over a field $k$, because it can be considered as a specialization of the diagonal carrying the motivic information at large. More precisely, let $k$ be an algebraically closed field, let $d$ be the dimension of $X$, and let $K=k(X)$ be the function field on $X$. Consider a pull-back homomorphism
  $$
  \Phi :CH^d(X\times X)\to CH^d(X_K)
  $$
induced by the embedding of the generic point $\eta =\Spec (K)$ into $X$. The kernel of $\Phi $ is generated by correspondences supported on $Z\times X$, where $Z$ runs Zariski closed subschemes in $X$ different from $X$ itself, see \cite{KMP}. Hence, various motivic effects, given originally in terms of correspondences, i.e. cycle classes in $CH^d(X\times X)$, can be expressed in terms of zero-cycles on $X_K$, modulo motives of varieties of dimension $<d$.

Assume, for example, that $X$ is a surface of general type over an algebraically closed field $k$, and the second Weil cohomology group $H^2(X)$ is algebraic. Let $\Delta _X$ be the diagonal on $X\times X$. Its specialization
  $$
  P_{\eta }=\Phi (\Delta )
  $$
is the generic point $\eta $ viewed as a zero-cycle on $X_K$. Fix now a $k$-rational point $P_0$ on $X$. Let $\Omega $ be a universal domain containing $k$ and embed $K$ into $\Omega $ over $k$. In the paper we will show, see Corollary \ref{Bloch}, that if $P_{\eta }$ is rationally equivalent to $P_0$ on $X_{\Omega }$ then any point $P$ is rationally equivalent to any other point $Q$ on $X_{\Omega }$, i.e. Bloch's conjecture hold's for $X_{\Omega }$. As Bloch's conjecture is equivalent to finite-dimensionality of the motive $M(X_{\Omega })$, we see that the above specialization map $\Phi $ allows to reformulated motivic effects at large in terms of rational equivalence between two concrete points on $X_{\Omega }$.

Certainly, it is still not easy to prove (or disprove) rational equivalence between the above points $P_{\eta }$ and $P_0$. One of the problems here consists of the lack of rational curves on surfaces of general type with algebraic $H^2(X)$. However, we believe that any further progress towards Bloch's conjecture must involve analysis of a possibility of an explicit rational deformation of $P_{\eta }$ into $P_0$ on the surface $X_{\Omega }$.

\medskip

The above picture can now be generalized as follows. Let $X$ be a smooth projective variety of dimension $d$ over an algebraically closed field $k$. To any zero-cycle $Z=\sum _in_iP_i$ on $X$ one can define its transcendence degree as the maximum of transcendence degrees of the residue fields $k(P_i)$. The transcendence degree of a zero-cycle class $\alpha \in CH^d(X)$ is the exact lower bound of the transcendence degrees of representatives of $\alpha $. Then the motive $M(X)$ is a direct summand of motives of varieties of dimensions $<d$, twisted by Lefschetz motives, if and only if the transcendence degree of any zero-cycle class $\alpha \in CH^d(X)$ is strictly smaller than $d$.

A nice thing is that the last assertion is also equivalent to the fact that there exists a point $P$ of transcendence degree $d$ on $X_{\Omega }$, rationally equivalent to a zero-cycle on $X$ whose transcendence degree is strictly smaller than $d$. More precisely, we prove the following theorem (see Theorem \ref{theor-maintrdeg} in the text below):

\medskip

\begin{itemize}

\item[]{}

{\it For any smooth projective variety $X$ of dimension $d$ over $k$ the following conditions are equivalent:

\begin{enumerate}

\item[]{(i)}
the class of the diagonal $\Delta _X$ is balanced;

\item[]{(ii)}
the Chow motive of $X$ is a direct summand of a sum of motives of varieties of dimension strictly smaller than $d$;

\item[]{(iii)}
the transcendence degree of any zero-cycle class on $X_{\Omega }$ is strictly less than $d$;

\item[]{(iv)}
there exists a closed point on $X_{\Omega }$ whose transcendence degree is $d$ but the transcendence degree of its class modulo rational equivalence is strictly less than $d$.

\end{enumerate}

}

\end{itemize}

\medskip

\bigskip

{\sc Acknowledgements.} The first author was partially supported by the grants RFBR 08-01-00095, NSh-4713.2010.1 and MK-297.2009.1. Both authors are thankful to Artiom Brazovsky for the hospitality in his country-house in Zadomlya (Belarus) where the draft version of this paper has been designed in August 2010.

\bigskip

\bigskip

\section{Some motivic lemma}

Below we will use the notation from \cite{GG3}. In particular, all Chow groups will be with coefficients in $\QQ $, unless the cases when integral coefficients will be subscripted by $\ZZ $. The category of Chow motives $\CHM $ over a field $k$ will be contravariant. That is, if $X$ and $Y$ are two smooth projective varieties over $k$, and $X$ is decomposed in to its connected components $X_j$, then the group of correspondences
  $$
  \Corr ^m(X,Y)
  $$
of degree $m$ from $X$ to $Y$ is a direct sum of the Chow groups
  $$
  CH^{e_j+m}(X_j\times Y)\; ,
  $$
where $e_j$ is the dimension of $X_j$. The composition of two correspondences $f\in \Corr ^n(X,Y)$ and $g\in \Corr ^m(Y,Z)$ is standard
  $$
  g\circ f={p_{13}}_*(p_{12}^*(f)\cdot p_{23}^*(g))\; ,
  $$
where the dot stands for the intersection of cycle classes in the sense of \cite{Fulton}. We also have a contravariant functor $M$ from the category of smooth projective varieties over $k$ to the category $\CHM $ sending any variety $X$ to its motive
  $$
  M(X)=(X,\Delta _X,0)\; ,
  $$
where $\Delta _X$ is the diagonal class of $X$, and any morphism $f:X\to Y$ maps to the class of transposition of its graph
  $$
  \Gamma _f^{\rm t}\in \Corr ^0(X,Y)\; .
  $$
The category of Chow motives $\CHM $ is tensor, with the tensor product induced by the products of varieties. The unite motive
  $$
  \uno =(\Spec (k),\Delta _{\Spec (k)},0)
  $$
and the Lefschetz motive
  $$
  \Le =(\Spec (k),\Delta _{\Spec (k)},-1)
  $$
are related by the formula
  $$
  M(\PR ^1)=\uno \oplus \Le \; .
  $$
For any positive integer $m$ let $\Le ^m$ be the $m$-fold tensor power of the Lefschetz motive $\Le $, let $\Le ^0=\uno $ and let $\Le ^{-m}=(\Le ^{-1})^{\otimes {-m}}$, where
  $$
  \Le ^{-1}=(\Spec (k),\Delta _{\Spec (k)},-1)\; .
  $$

Further details on Chow motives can be found, for example, in \cite{GG3}.

\bigskip

The next notion we need is the notion of balancing. Let $X$ and $Y$ be two equi-dimensional varieties over $k$. Similarly to \cite{Barbieri Viale}, we say that a correspondence $\alpha\in CH^m(X,Y)$ is {\it balanced on the left} (respectively, {\it on the right}) if there exists an equi-dimensional Zariski closed subscheme $Z\subset X$, such that
  $$
  \dim (Z)<\dim (X)\; ,
  $$
and an algebraic cycle $\Gamma $ on $X\times Y$, such that $[\Gamma ]=\alpha $ in $CH^m(X,Y)$ and the support of $\Gamma $ is contained in $Z\times X$ (respectively, in $X\times Z$). The subscheme $Z$ will be called a {\it pan} of balancing. We say that $\alpha \in CH^m(X,Y)$ is {\it balanced} if $\alpha =\alpha _1+\alpha _2$, where $\alpha _1$ is balanced on the left, and $\alpha_2$ is balanced on the right.

\bigskip

Balancing was discovered in \cite{Bloch} and \cite{BlochSrinivas}. It is a motivic notion and can be restated in purely motivic terms:

\bigskip

\begin{lemma}
\label{lemma-suppleft}
Let $X$ and $Y$ be equidimensional smooth projective varieties over $k$, and let $\alpha \in CH^m(X,Y)$. Then $\alpha $ is balanced on the left if and only if there exists an equidimensional smooth projective variety $Z$ over $k$ with
  $$
  \dim (Z)<\dim (X)\; ,
  $$
such that $\alpha $ factors through $M(Z)$, that is $\alpha$ is a composition
  $$
  M(X)\lra M(Z)\lra M(Y)\otimes \Le ^{-m}\; .
  $$
Symmetrically, the correspondence $\alpha $ is balanced on the right if and only if there exists an equidimensional smooth projective variety $Z$ over $k$ with
  $$
  n=\dim(Y)-\dim(Z)>0\; ,
  $$
such that $\alpha $ is a composition
  $$
  M(X)\lra M(Z)\otimes \Le ^{n-m}\lra M(Y)\otimes \Le^{-m}\; .
  $$
\end{lemma}

\begin{pf}
If $m=0$ and the closed subscheme $Z$ is smooth, then the lemma is just obvious. Indeed, let $i:Z\hra X$ be the closed embedding, and let $\Gamma _i^{\rm t}$ be the transpose of the graph of the embedding $i$. If $\alpha $ is balanced on the left then it can be considered as a correspondence of degree zero from $Z$ to $Y$. Therefore, the correspondence $\alpha $ from $X$ to $Y$ is a composition of the correspondence $\Gamma _i^{\rm t}$ with $\alpha $ as a correspondence from $Z$ to $Y$.

The detailed proof of the lemma when $Z$ is not necessarily smooth and $m\neq 0$ is given in \cite{GG3}.
\end{pf}

In the next section we will introduce the transcendence degree of a zero cycle on a smooth projective variety and we will show how it is related to balancing of the diagonal, and so the above motivic factorizations from Lemma \ref{lemma-suppleft}.

\section{Transcendence degree of zero-cycles}

First we need to recall some well-known things from the theory of schemes.

\medskip

Let $k$ be a field, and let $X$ be an algebraic scheme over $k$. Let $k\subset K$ be a field extension. Recall that a $K$-point on $X$ is a morphism of schemes $P:\Spec (K)\to X$ over $\Spec (k)$. A subextension $k\subset L\subset K$ is a field of definition of the point $P$ if there exists a morphism
  $$
  p_L : \Spec (L)\lra X\; ,
  $$
such that the following diagram is commutative
  $$
  \xymatrix{
  \Spec (K) \ar[ddrr]^-{} \ar[rrrr]^-{P}  & & & & X  \\ \\
  & & \Spec (L) \ar[rruu]_-{P_L} & &
  }
  $$
as a diagram over $\Spec (k)$.

Let $\xi _P$ be the image of the unique point in $\Spec (L)$ with respect to the morphism $P$, and let $k(\xi _P)$ be the residue field of the point $\xi _P$ on the scheme $X$. Then $k(\xi _P)$ is the minimal field of definition of the point $P$, i.e. the initial object in the category of fields of definition of the point $P$, because $k(\xi _P)\hra L$ and the above morphism $P_L$ factors through the natural morphism $\Spec (k(\xi _P))\to X$.

By definition, the transcendence degree of the point $P$ over the ground field $k$ is the transcendence degree of the field $k(\xi _P)$ over $k$:
  $$
  \trdeg (P/k)=\trdeg (k(\xi _P)/k)\; .
  $$

Thus, the transcendence degree $\trdeg (P/k)$ is the transcendence degree of the minimal field of definition of the point $P$ over the ground field $k$.

Notice that if $k\subset L\subset K$ is a field subextension then one has a commutative diagram
  $$
  \xymatrix{
  \Spec (K) \ar[ddrr]^-{} \ar[rrrr]^-{Q_L}  & & & &
  X \times _{\Spec (k)}\Spec (L) \ar[lldd]_-{} \\ \\
  & & \Spec (L)  & &
  }
  $$
Notice that the transcendence degree $\trdeg (Q_L/L)$ can be different from the transcendence degree $\trdeg (P/k)$. For example, $\trdeg (P_{k(\xi _P)}/k(\xi _P))=0$.

\medskip

Let $Y$ be the Zariski closure of the schematic point $\xi _P$ in $X$. Then $Y$ is a closed irreducible subscheme in $X$ and
  $$
  \trdeg (P/k)=\dim _k(Y)\; .
  $$
It follows, in particular, that
  $$
  \trdeg (P/k)\leq \dim _k(X)\; .
  $$

\medskip

Now we are going to introduce the notion of a transcendence degree of a zero-cycle on a variety. Let $\Omega $ be a universal domain containing $k$. Suppose $X$ is an equidimensional variety, and let $d$ be the dimension of $X$.

\medskip

\begin{definition}
A {\rm transcendence degree} $\trdeg (\alpha /k)$ of a zero-cycle class $\alpha \in CH^d(X_{\Omega })$ over $k$ is the minimal natural number $n$, such that there exists a zero-cycle
  $$
  Z=\sum _in_iP_i
  $$
on $X_{\Omega }$ representing the class $\alpha $ with the property
  $$
  \trdeg (P_i/k)\leq n
  $$
for all $i$.
\end{definition}

\medskip

The following properties of the transcendence degree for zero-cycles follow directly from the above definition.

\begin{lemma}
\label{lemma-propertiestrdeg}
Let $X$ be an equidimensional variety over $k$ of dimension $d$. Then the following is true:

  \begin{enumerate}

  \item[]{(i)}
  for any element $\alpha\in CH^d(X_{\Omega })$ one has
    $$
    \trdeg (\alpha /k)\le d;
    $$

  \item[]{(ii)}
  for all elements $\alpha ,\beta \in CH^d(X_{\Omega })$ we have that
    $$
    \trdeg((\alpha+\beta)/k)\le \max\{\trdeg(\alpha/k)\, ,\; \trdeg(\beta/k)\};
    $$

  \item[]{(iii)}
  given a field subextension $k\subset K\subset \Omega $ and an element $\beta \in CH^d(X_K)$, we have an inequality
    $$
    \trdeg(\beta\,_{\Omega }/k)\le \trdeg(K/k)\; .
    $$

  \end{enumerate}

\end{lemma}

\begin{remark}
{\rm Not any cycle class $\alpha \in CH^d(X_{\Omega })$ is equal to $\beta\, _{\Omega }$, for some $\beta \in CH^d(X_K)$ and $K$ with $\trdeg (K/k)=\trdeg (\alpha /k)$. Let, for example, $X$ be a smooth projective curve of genus at least two. Then there exists a point $P$ of transcendence degree at least two on the Jacobian variety $\Jac (X)$ of $X$ over $k$. Let $\alpha $ be a cycle class in the Chow group $CH^1(X_{\Omega })_0$ of degree zero $0$-cycles on the curve $X$ corresponding to the point $P$ under the isomorphism
  $$
  CH^1(X_{\Omega })_0=\Jac (X)_{\Omega }\; .
  $$
Then $\trdeg (\alpha )\leq 1$ because $\dim (X)=1$. Suppose now that $\alpha $ comes from an element $\beta \in CH^1(X_K)_0$ by means of the scalar extension from $K$ to $\Omega $, where $\trdeg (K/k)=1$. Since the isomorphism between the Chow group of degree zero $0$-cycles and the Jacobian commutes with scalar extensions of the ground field, the point $P$ must be defined over $K$, which is impossible as $\trdeg (P/k)=2$.
}
\end{remark}

\medskip

We will also use the following fact.

\begin{lemma}\label{lemma-motivictrdeg}
Let $X$ and $Y$ be two smooth projective equidimensional varieties over $k$, let $d=\dim (X)$, $e=\dim (Y)$ and assume $e<d$. Let $\varphi $ be a correspondence of degree $d-e$ from $Y$ to $X$, that is $\varphi$ is a morphism of Chow motives
  $$
  M(Y)\otimes \Le ^{\otimes (d-e)}\lra M(X)\; .
  $$
Then for any element $\alpha \in CH^e(Y_{\Omega })$ one has
  $$
  \trdeg ((\varphi _{\Omega })_*(\alpha )/k)\leq \trdeg (\alpha /k)\; .
  $$
\end{lemma}

\begin{pf}
Let $\sum _in_iP_i$ be a zero-cycle on $Y_{\Omega }$, such that
  $$
  \alpha =\left[\sum _in_iP_i\right]\; ,
  $$
and
   $$
   \trdeg (P_i/k)\le \trdeg (\alpha /k)
   $$
for all $i$. By linearity of the push-forward homomorphism $(\varphi _{\Omega })_*$ and also by Lemma~\ref{lemma-propertiestrdeg} (ii), it is enough to show that
  $$
  \trdeg ((\varphi _{\Omega })_*[P_i]/k)\leq \trdeg (P_i/k)
  $$
holds true for all indices $i$.

Let $P$ be one of the points $P_i$. By the definition of a transcendence degree of a point there exists a field $K$, such that
  $$
  \trdeg (K/k)=\trdeg (P/k)\; ,
  $$
and a point $W\in Y(K)$, such that
  $$
  W_{\Omega }=P\; .
  $$
Moreover,
  $$
  (\varphi _{\Omega })_*[P]=((\varphi _K)_*[W])_{\Omega }\; .
  $$
Since $(\varphi _K)_*[W]\in CH^d(X_K)$, by Lemma~\ref{lemma-propertiestrdeg}(3), we see that
  $$
  \trdeg ((\varphi _{\Omega })_*[P]/k)\leq \trdeg (K/k)=\trdeg (P/k)\; ,
  $$
which completes the proof.
\end{pf}

\begin{remark}
{\rm Certainly, one can also define the notion of a transcendence degree for all closed irreducible subschemes in $X_{\Omega }$ and, respectively, for elements in Chow groups $CH^p(X_{\CC})$ of arbitrary codimension $p$. Moreover, analogs of Lemma~\ref{lemma-propertiestrdeg} (ii) (iii) and Lemma~\ref{lemma-motivictrdeg} imply that a transcendence degree is also well-defined for elements in Chow groups of Chow motives over $k$, and that this transcendence degree does not increase under taking push-forwards with respect to morphisms between Chow motives over $k$.
}
\end{remark}

\bigskip

Now we are ready to prove our main statement.

\medskip

\begin{theorem}
\label{theor-maintrdeg}
Let $X$ be an irreducible smooth projective variety over $k$ of dimension $d$. The following conditions are equivalent:

\begin{enumerate}

\item[]{(i)}
the class of the diagonal $\Delta _X$ is balanced in $CH^d(X\times X)$;

\item[]{(ii)}
the Chow motive $M(X)$ is isomorphic to a direct summand of the motive
  $$
  M(Y_1)\oplus (M(Y_2)\otimes\Le^{d-e})\; ,
  $$
where $Y_1$ and $Y_2$ are equidimensional smooth projective varieties over $k$ whose dimensions are strictly less than $d$, and $e$ is the dimension of the variety $Y_2$;

\item[]{(iii)}
any element $\alpha \in CH^d(X_{\Omega })$ satisfies
  $$
  \trdeg (\alpha /k)<d\; ;
  $$

\item[]{(iv)}
there exists a closed point $P\in X_{\Omega }$, such that $\trdeg (P/k)=d$ and
  $$
  \trdeg ([P])<d\; ,
  $$
where $[P]$ is the class of the point $P$ in $CH^d(X_{\Omega })$.

\end{enumerate}

\end{theorem}

\begin{pf}

\bigskip

\noindent $(i)\Rightarrow(ii)$

\medskip

Suppose that $[\Delta _X]=\alpha _1+\alpha _2$, where $\alpha _1$ is balanced on the left and $\alpha _2$ is balanced on the right. By Lemma~\ref{lemma-suppleft}, there exist two equidimensional varieties $Y_1$ and $Y_2$ as in (2), and factorizations of $\alpha _1$ and $\alpha _2$, so that $\alpha $ factorizes like this:
  $$
  M(X)\to M(Y_1)\to M(X),\quad M(X)\to M(Y_2)\otimes \Le ^{d-e}\to M(X)\; .
  $$
Put $M:=M(Y_1)\oplus (M(Y_2)\otimes\Le^{d-e})$. Then the identity morphism from $M(X)$ to itself factors through $M$, thus, $M(X)$ is a direct summand in $M$.

\bigskip

\noindent $(ii)\Rightarrow(iii)$

\medskip

Looking at the Chow groups of the motives involved in the decomposition
  $$
  M(X)\lra M(Y_1)\oplus (M(Y_2)\otimes\Le^{d-e})\lra M(X)
  $$
we see that all elements in $CH^d(X_{\Omega })$ are push-forwards with respect to the morphism
  $$
  M(Y_2)\otimes\Le^{d-e}\lra M(X)\; ,
  $$
as
  $$
  CH^d(M((Y_1)_{\Omega }))=CH^d((Y_1)_{\Omega })=0
  $$
because $e<d$. Then (iii) follows from Lemma~\ref{lemma-motivictrdeg} and Lemma~\ref{lemma-propertiestrdeg} (i).

\bigskip

$(iii)\Rightarrow(iv)$

\medskip

This is just obvious.

\bigskip

$(iv)\Rightarrow(i)$

\medskip

Let $\sum _in_iP_i$ be a zero-cycle on $X_{\Omega }$, such that
  $$
  [P]=\left[\sum _in_iP_i\right]\; ,
  $$
and
  $$
  \trdeg (P_i/k)<d\; .
  $$
By definition of a transcendence degree, there are field extensions
  $$
  K\subset \Omega \quad \hbox{and}\quad K_i\subset \Omega
  $$
over $k$, and points
  $$
  W\in X(K)\; ,\; \; W_i\in X(K_i)\; ,
  $$
such that
  $$
  W_{\Omega }=P\; ,\; \; \; (W_i)_{\Omega }=P_i\; ,
  $$
the fields $K$ and $K_i$ are finitely generated over $k$ with
  $$
  \trdeg (K/k)=d\quad \hbox{and}\quad  \trdeg(K_i/k)<d\; .
  $$

Let $L$ be the composite of the fields $K$ and $K_i$ in $\Omega $. As
  $$
  [P]=\left[\sum _in_iP_i\right]\in CH^d(X_{\Omega })
  $$
and all involved Chow groups are with coefficients in $\QQ $, one has
  $$
  [W_L]=\left[\sum _in_i(W_i)_L\right]\in CH^d(X_L)\; ,
  $$
see \cite{Bloch}, page 1.21.

Let now $V$ be a smooth irreducible quasi-projective variety over $k$, such that
  $$
  k(V)=K\; .
  $$
Then we also have a rational dominant morphism
   $$
   f:V{\dashrightarrow }X\; ,
   $$
which coincides at the generic point with the morphism $W:\Spec(K)\to X$.

Similarly, for each $i$, we have a smooth irreducible quasi-projective variety $V_i$ with $k(V_i)=K_i$, and a rational dominant morphism
  $$
  f_i:V_i{\dashrightarrow }X
  $$
inducing the morphism $W_i:\Spec(K_i)\to X$ at the generic point.

Shrinking the varieties $V$ and $V_i$ to Zariski open subsets one can think that the above morphisms $f$ and $f_i$ are all regular.

We also need a smooth irreducible quasi-projective variety $Z$ over $k$ with dominant regular morphisms $g:Z\to V$ and $g_i:Z\to V_i$, such that the function field $k(Z)$ coincides with $L$.

For any regular morphism $h$ let $\Gamma _h$ be the graph of $h$. Shrinking $Z$ to a non-empty Zariski open subset if necessary, we have that
  $$
  [\Gamma _{fg}]=\left[\sum _in_i\Gamma _{f_ig_i}\right]
  $$
in the group $CH^d(Z\times X)$, because the analogous rational equivalence holds over the generic point of $Z$, which is $\Spec (L)$, see above.

Notice that
  $$
  \dim (V)=\trdeg (K/k)=d\; .
  $$

Let $T\subset Z$ be a generic $d$-dimensional multiple hyperplane section of $Z$. The scheme $T$ is irreducible by Bertini's theorem, and the restrictions
  $$
  h:=g|_T:T\to W,\quad h_i:=(g_i)|\, _T:T\to W_i
  $$
are still dominant. By taking pull-backs in Chow groups with respect to the embedding $T\times X\to Z\times X$, we obtain
  $$
  [\Gamma _{fh}]=\left[\sum _in_i\Gamma _{f_ih_i}\right]
  $$
in the group $CH^d(T\times X)$.

Since $\dim (T)=d$ and the composition $fh:T\to X$ is dominant, we see that $fh$ is generically finite. Thus, shrinking $T$ to a non-empty open subset, we may assume that the morphism $fh$ is a finite surjective morphism from $T$ onto a non-empty open subset $U$ in $X$.

Now we use push-forwards in Chow groups with respect to the finite morphism
  $$
  fh\times \id_X:T\times X\lra U\times X\; .
  $$
From the above equality we obtain that
  $$
  (fh\times \id _X)_*[\Gamma _{fh}]=
  (fh\times \id _X)_*\left[\sum _in_i\Gamma _{f_ih_i}\right]
  $$
in the group $CH^d(U\times X)$.

Set-theoretically,
  $$
  (fh\times \id _X)(\Gamma _{fh})=\Delta _X\cap (U\times X)\; .
  $$
The closure of $(fh\times \id _X)(\Gamma _{f_ih_i})$ in $X\times X$ is contained in $X\times \overline {f_i(W_i)}$, where $\overline {f_i(W_i)}$ is the Zariski closure of $f_i(W_i)$ in $X$.

Since
  $$
  \dim (V_i)=\trdeg (K_i/k)<d
  $$
and all the Chow groups are with rational coefficients, we see that $\Delta _X$ is balanced.
\end{pf}

\medskip

\noindent {\bf Remark.} The equivalence $(i)\Leftrightarrow(ii)$ was actually proved in \cite{GP2} but we included it in the theorem for the convenience of the reader.

\medskip

\section{An example}

An important thing in Theorem \ref{theor-maintrdeg} is that (iv) implies (i). Let us illustrate this by an example.

Let $X$ be a smooth projective surface over $\CC $, of general type and with $p_g=0$. Recall, that Bloch's conjecture predicts that for any two closed points $P$ and $Q$ on $X_{\CC }$ the point $P$ is rationally equivalent to $Q$. This conjecture is a codimension $2$ case of the Bloch-Beilinson paradigma for algebraic cycles, and it is highly inaccessible. It is known for surfaces with the Kodaira dimension $<2$, \cite{BKL}, for finite quotients of products of curves, \cite{Kimura}, and for surfaces of general type (which are not finite quotients of products of curves) in \cite{InoseMizukami}, \cite{Barlow} and \cite{Voisin}.

Let now $k$ be the algebraic closure in $\CC $ of the minimal field of definition of the surface $X$, and let $K=k(X)$ be the function field of $X$ over $k$. Let $\eta =\Spec (K)$ be the generic point of $X$, and let $P_{\eta }$ be the corresponding $K$-rational closed point on $X_K$. Theorem \ref{theor-maintrdeg} implies the following corollary:

\medskip

\begin{corollary}
\label{Bloch}
Bloch's conjecture holds for $X$ if and only if there exist an embedding of $K$ into $\CC $ over $k$, and a $k$-rational point $P$ on $X$, such that the above closed $K$-rational point $P_{\eta }$ is rationally equivalent to $P$ on $X$ as a variety over $\CC $.
\end{corollary}

\medskip

This can be made absolutely explicit in the case of Godeaux surfaces, for which Bloch's conjecture was proved by C.Voisin in \cite{Voisin}. Namely, let $\mu _5$ be the group of $5$-th roots of the unit in $\CC $, and let $\epsilon $ be a primitive root in it. The group $\mu _5$ acts on $\PR ^3$ by the rule:
  $$
  [x_0:x_1:x_2:x_3]\mapsto [x_0:\epsilon x_1:\epsilon ^2x_2:\epsilon ^3x_3]
  $$
Let $f=f(x)$ be a $\mu _5$-invariant smooth quintic form in $\PR ^3$, and let $Y=Z(f)$ be the set of zeros of $f$ in $\PR ^3$. Since $f$ is $\mu _5$-invariant, the group $\mu _5$ acts on $Y$. Assume, in addition, that $Y$ does not contain the four fixed points of the action of $\mu _5$ on $\PR ^3$. Then the quotient surface
  $$
  X=Y/\mu _5
  $$
is non-singular, and it is called a Godeaux surface. It is well known that $p_g=q=0$ for such $X$, see \cite{Reid}.

Take now two transcendental complex numbers which are algebraically independent over $\QQ $, say $e$ and $e^{\pi }$, see \cite{Nesterenko}. Let $\alpha $ be one of the zeros of the polynomial obtained by substitution of the coordinates $e$ and $e^{\pi }$ in to the affinized form $f$. Then $P_{\eta }$ can be represented as the class of the point
  $$
  (e,e^{\pi },\alpha )\in \CC ^3
  $$
under the quotient-map $Y\to X$.

Then Voisin's result says that the point $P_{\eta }$ is rationally equivalent to a point in $X(\bar \QQ )$. The specificity of Corollary \ref{Bloch} is that it says that the above rational equivalence between two single points on $X(\CC )$ is the only reason for vanishing of the whole Albanese kernel in this situation.

We believe that this observation can be useful in approaching to Bloch's conjecture in some concrete contexts, such as Mumford's fake surface, see \cite{Mumford}. Recall that such surfaces were recently classified in \cite{PrasadYeung}.

\bigskip

WARNING. It would be a temptation to find a rational curve through the points $P_{\eta }$ and $P_0$ on the Godeaux surface $X$ over $\CC $. The first problem is that $X$ is a surface of general type whose discrete invariants vanish, so that one can expect only a few rational curves on $X_{\CC }$. But this is not yet the main trouble. The main difficulty is that no rational curves can pass through $P_{\eta }$ at all.

\medskip

Indeed, let $X$ be a smooth projective surface over the ground subfield $k$ in $\Omega $. Let $P_{\eta }$ be a closed point of transcendence degree $2$ on $X_{\Omega }$. Suppose there exists a field subextension $k\subset K\subset \Omega $, a point $P:\Spec (K)\to X$ with $\trdeg (P/k)=2$, and a rational curve $C$ on $X_K$ passing through the point $P$. Let us show that $X$ is uniruled then. Without loss of generality one can assume that $K$ is finitely generated over $k$. Let $Y$ be an irreducible variety over $k$, such that $K$ is the function field of $Y$ over $k$.The rational curve $C\subset X_K$ induces a morphism
  $$
  \phi :\PR ^1_K\lra X_K
  $$
which induces a rational morphism
  $$
  \PR ^1\times _kY\dashrightarrow X\times _kY\; ,
  $$
such that
  $$
  f\times _Y\Spec (K)=\phi \; .
  $$
The point $p$ gives a morphism
  $$
  \Spec (K)\to \PR ^1\times _kK
  $$
over $K$. This corresponds to some rational section of the projection
  $$
  \PR ^1\times Y\to Y\; .
  $$
The morphism
  $$
  \Spec (K)\to \PR ^1_K\to X_K\to X
  $$
sends the unique point in $\Spec (K)$ in to the generic point of $X$ because $\trdeg (p/k)=2$. Therefore, the composition
  $$
  Y\dashrightarrow \PR ^1\times _kY\stackrel{f}{\dashrightarrow }X\times _kY
  \stackrel{p_X}{\dashrightarrow }X
  $$
is dominant, where $p_X$ is the projection onto $X$. It follows that the morphism
  $$
  \PR ^1\times _kY\lra X
  $$
is dominant as well. Moreover, the induced map
  $$
  \PR ^1_K\lra X_K
  $$
gives a birational isomorphism with its image. It follows that this image is a curve in $X_K$. Hence, the map
  $$
  \PR ^1\times _kY\lra X
  $$
does not factor through the projection $\PR ^1\times Y\to Y$. Hence, at least for one point $y\in Y$ the induced map
  $$
  \PR ^1_y\lra X
  $$
is not constant. Hence, $X$ is uniruled by \cite[1.3.4]{Kollar}

\medskip

Thus, if we could have a rational curve through $P_{\eta }$ on a smooth projective surface $X_{\CC }$, of general type with $p_g=0$, then immediately we would get a contradiction as such a surface is very far from to be uniruled.

\medskip

This shows that in order to find a precise rational equivalence between $P_{\eta }$ and $P_0$ we need to find more than one curves of genus $>0$ on the Godeaux surface $X$, and rational functions of them, which will provide a suitable zero-poles cancelation for their principle divisors.

\bigskip

\begin{small}

\end{small}

\bigskip

\bigskip

\begin{small}

{\sc Steklov Mathematical Institute, Gubkina str. 8,
119991, Moscow, Russia}

\end{small}

\medskip

\begin{footnotesize}

{\it E-mail address}: {\tt gorchins@mi.ras.ru}

\end{footnotesize}

\bigskip

\begin{small}

{\sc Department of Mathematical Sciences, University of Liverpool,
Peach Street, Liverpool L69 7ZL, England, UK}

\end{small}

\medskip

\begin{footnotesize}

{\it E-mail address}: {\tt vladimir.guletskii@liverpool.ac.uk}

\end{footnotesize}


\begin{thebibliography}{9999}

\bibitem{Barbieri Viale}
L.Barbieri Viale. Balanced varieties. Algebraic $K$-theory and its applications (Trieste, 1997), 298 - 312, World Sci. Publ., River Edge, NJ, 1999

\bibitem{Barlow}
R. Barlow. Rational equivalence of zero cycles for some more surfaces with $p_g=0$. Inventiones Mathematicae 79 (1985) 303 - 308

\bibitem{Bloch}
S. Bloch. Lectures on algebraic cycles. Duke Univ. Math. Series IV, 1980

\bibitem{BKL}
S. Bloch, A. Kas, D. Lieberman. Zero cycles on surfaces with $p_g=0$. Compositio Math. 33 (1976), 135 - 145

\bibitem{BlochSrinivas}
S. Bloch, V. Srinivas. Remarks on correspondences and algebraic cycles. Amer. J. Math. 105 (1983) 1235 - 1253

\bibitem{Fulton}
W. Fulton. Intersection theory. Ergebnisse der Mathematik und ihrer Grenzgebiete, Vol. 3, No. 2, Springer-Verlag (1984)


\bibitem{GG3}
S.Gorchinskiy, V.Guletski\u \i . Non-trivial elements in the Abel-Jacobi kernels on higher dimensional varieties. Preprint 2010 (see ArXiv)

\bibitem{GP2}
V. Guletski\u {\i }, C. Pedrini. Finite-dimensional motives and the conjectures of Beilinson and Murre. $K$-Theory, Vol. 30, No. 3 (2003), 243
- 263

\bibitem{InoseMizukami}
H. Inose, M. Mizukami. Rational equivalence of zero-cycles on some surfaces with $p_g=0$. Math. Ann. 244 (1979) 205 - 217

\bibitem{KMP}
B.Kahn, J.Murre, C.Pedrini. On the transcendental part of the motive of a surface. Algebraic cycles and motives. Volume 2. London Mathematical Society. Lecture Note Series 344. Cambridge University Press 2007

\bibitem{Kimura}
S.-I. Kimura. Chow groups are finite dimensional, in some sense. Math. Ann. Vol. 331, no. 1 (2005) 173 -201

\bibitem{Kollar}
J.Koll\'ar. Rational curves on algebraic varieties. Ergebnisse der Mathematik und ihrer Grenzgebiete. 3. Folge. A Series of Modern Surveys in Mathematics 32. Springer-Verlag, Berlin, 1996

\bibitem{Mumford}
D. Mumford. An algebraic surface with $K$ ample, $K^2=9$, $p_g=q=0$. American Journal of Mathematics 101 (1) (1979) 233 - 244

\bibitem{Nesterenko}
Y.Nesterenko. Modular Functions and Transcendence Problems. Comptes rendus de l'Acad\'emie des sciences S\'erie 1 322 (10) (1996) 909 - 914

\bibitem{PrasadYeung}
G. Prasad, S. Yeung. Fake projective planes. Inventiones Mathematicae 168 (2) (2007) 321 - 370

\bibitem{Reid}
M.Reid. Campedelli versus Godeaux. Problems in the theory of surfaces and their classification (Cortona, 1988) 309 - 365. Sympos. Math., XXXII, Academic Press, London, 1991

\bibitem{Voisin}
C. Voisin. Sur les z\'ero-cycles de certaine hypersurfaces munies d'un automorphisme. Ann. Scuola Norm. Sup. Pisa Ck. Si. (4) 19 (1992) 473 - 492

\end{thebibliography}
\end{document}